\documentclass[11pt]{article}
\usepackage{amsmath,amssymb,amsthm,amsfonts,epsfig}
\input{xy}                                                                      
\xyoption{all}
\CompileMatrices 

\newcommand{\com}[1]{
}

\newcommand{\us}[1]{{\upshape{#1}}}
\newcommand{\bp}{\begin{pmatrix}}
\newcommand{\ep}{\end{pmatrix}}
\newcommand{\be}{\begin{equation}}
\newcommand{\ee}{\end{equation}}
\newcommand{\bs}{\begin{split}}
\newcommand{\es}{\end{split}}
\newcommand{\bc}{\begin{center}}
\newcommand{\ec}{\end{center}}
\newcommand{\ed}{{\rm d}}
\newcommand{\w}{{\mathchoice{\,{\scriptstyle\wedge}\,}{{\scriptstyle\wedge}}    
      {{\scriptscriptstyle\wedge}}{{\scriptscriptstyle\wedge}}}}                
  
\newcommand{\ol}{\overline}
\newcommand{\liealgebra}[1]{{\mathfrak {#1}}}

\newcommand{\spin}{\liealgebra{spin}}
\newcommand{\liegroup}[1]{{\operatorname{#1}}}

\newcommand{\SO}{\liegroup{SO}}
\newcommand{\SU}{\liegroup{SU}}
\newcommand{\Un}{\liegroup{U}}

\newcommand{\Spin}{\liegroup{Spin}}
\DeclareMathOperator{\stab}{Stab}
\renewcommand{\Re}{\operatorname{Re}}
\renewcommand{\Im}{\operatorname{Im}}
\newcommand{\om}{\omega}
\newcommand{\Om}{\Omega}

\newcommand{\gro}{{\rm G}(2,\Oc)}
\newcommand{\gris}{{\rm Gr}^{0}(2,\C^4)}

\newcommand{\R}{\mathbb R}
\newcommand{\C}{\mathbb C}

\newcommand{\Z}{\mathbb Z}
\newcommand{\Oc}{\mathbb O}
\newcommand{\s}[1]{{\mathbb S}^{#1}}

\newcommand{\mce}{\mathcal E}

\newcommand{\mch}{\mathcal H}
\newcommand{\mci}{\mathcal I}
\newcommand{\mcj}{\mathcal J}

\newcommand{\mcl}{\mathcal L}

\newcommand{\mco}{\mathcal O}

\newcommand{\mcv}{\mathcal V}
\newcommand{\mcw}{\mathcal W}

\newcommand{\ci}[1]{C^{\infty}({#1})}
\newcommand{\spn}{\Delta}

\newcommand{\ip}[2]{\langle {#1},\, {#2} \rangle}
\newcommand{\hra}{\hookrightarrow}
\newcommand{\p}{\varphi}
\newcommand{\G}{{\rm G}_2}

\newcommand{\ff}{{\rm I \negthinspace I}}

\newcommand{\cp}[1]{{\C \mathbb{P}^{#1}}}

\newcommand{\db}{\bar{\partial}}

\newcommand{\bx}{{\bf x}}
\newcommand{\bfp}{{\bf p}}

\newcommand{\ns}{\negthinspace}
\newcommand{\tr}{\rm tr}
\newcommand{\trp}{\; {}^t \ns}

\newcommand{\fb}{\bar{f}}
\newcommand{\hoz}{\mch'}

\newcommand{\vooz}{\mcv_1'}

\newcommand{\vtoz}{\mcv_2'}
\newcommand{\ezo}{\mce''}
\newcommand{\eoz}{\mce'}

\newcommand{\ffo}{ {\rm I}_1}
\newcommand{\ffe}{ {\rm I}_2}

\newcommand{\thod}{\theta_{\rm od}}
\newcommand{\thev}{\theta_{\rm ev}}
\usepackage[mathscr]{eucal}
\usepackage{bbm}
\usepackage{oldgerm}

\theoremstyle{plain}
\newtheorem{thm}{Theorem}[section]
\newtheorem{lem}[thm]{Lemma}
\newtheorem{cor}[thm]{Corollary}

\newtheorem{prop}[thm]{Proposition}

\newtheorem{rem}[thm]{Remark}
\theoremstyle{definition}
\newtheorem{defn}{Definition}[section]

\begin{document}
\title{Cayley Cones Ruled by $2$-Planes: Desingularization and Implications of the Twistor Fibration}
\author{Daniel Fox}
\date{\today}
\begin{abstract}
Cayley cones in the octonions $\Oc$ that are ruled by oriented 2-planes are equivalent to pseudoholomorphic curves in the Grassmannian of oriented $2$-planes $\gro$. The well known twistor fibration $\gro \to \s{6}$ is used to prove the existence of immersed higher-genus pseudoholomorphic curves in $\gro$.   Equivalently, this produces Cayley cones whose links are $\s{1}$-bundles over genus-$g$ Riemann surfaces.  When the degree of an immersed pseudoholomorphic curve is large enough, the corresponding $2$-ruled Cayley cone is the asymptotic cone of a non-conical $2$-ruled Cayley $4$-fold.
\end{abstract}

\thanks{This work grew out of a conversation with Dominic Joyce at the AIM Calibrations workshop in June 2006.    I would like to thank Ian McIntosh, Yinan Song, and especially Dominic Joyce for useful conversations while I worked on this paper.    I would like to thank the University of Oxford for its hospitality while I carried out this work.  This research was supported by NSF grant OISE-0502241.  
}

\maketitle

\tableofcontents


\begin{section}{Introduction}
We begin by showing that $r$-oriented $2$-ruled Cayley cones are equivalent to pseudoholomorphic curves in the Grassmannian of oriented $2$-planes in $\Oc$.  Results about such Cayley cones are then deduced by studying the geometry of the corresponding pseudoholomorphic curves.  Most significantly, the twistor fibration $\mcj:\gro \to \s{6}$ decomposes every pseudoholomorphic curve into an oriented branched minimal surface in $\s{6}$, a branched cover of Riemann surfaces, and a holomorphic line subbundle.  This decomposition allows an easy proof of the existence of immersed higher-genus pseudoholomorphic curves, and thus Cayley cones whose links have complicated topology.  By comparison, the proof of the existence of higher-genus embedded special Legendrian surfaces in $\s{5}$ by Haskins and Kapouleas \cite{hk} requires very involved and delicate analysis.  Neither the techniques nor the results of the current article seem to have any implications for special Legendrian surfaces in $\s{5}$.

The geometry of a pseudoholomorphic curve in $\gro$ defines a (possibly empty) class of deformations of the Cayley cone that preserve the properties of being Cayley and $2$-ruled but destroy the cone structure.  When the pseudoholomorphic curve has large negative degree then such deformations are abundant and will generically produce a smooth Cayley $4$-fold. 

The method for generating new pseudoholomorphic curves from old ones can also be viewed as a type of B\"acklund transformation in the sense of \cite{bgg}, \cite{clelland}, and \cite{ci}.  An article describing this is under preparation.  Similar transformations for certain harmonic maps have been described by Burstall in \cite{burst}. 
\end{section}

\begin{section}{$2$-ruled Cayley cones}\label{sec:streq}
Cayley submanifolds are a special class of minimal $4$-folds defined by first order partial differential equations.  We recall some of the basic facts and refer the reader to Appendix \ref{sec:spin7} and \cite{hl} for details\footnote{Beware that the choice of Cayley $4$-form in this article differs from \cite{hl} and most other conventions.}.  

Let $\Oc$ denote the octonions with standard coordinates $(x^1,\ldots,x^8):\Oc \to \R^8$.  There is a quadruple cross product on $\Oc$
\[
\begin{split}
&\Oc \times \Oc \times \Oc \times \Oc  \to \Oc\\
&(x,y,z,w) \mapsto x \times y \times z \times w 
\end{split}
\]
that is totally skew and whose components therefore define $4$-forms on $\Oc$.  Let $e^0_i=\frac{\partial}{\partial{x^i}}$, so that $\{ e^0_{i}\}$ provide an orthonormal basis of $\Oc$.  Assume that the octonionic structure is chosen so that $e^0_8=1 \in \Oc$ and the rest span $\Im(\Oc)$.  Define the $4$-forms 
\begin{align*}
\Phi(x,y,z,w)&=\ip{1}{x \times y \times z \times w}\\
\psi_m(x,y,z,w) &= \ip{e^0_m}{x \times y \times z \times w} \;\;{\rm for} \; m=1, \ldots, 7
\end{align*}
The stabilizer of $\Phi$ is $\Spin(7) \subset \SO(8)$ and so $\Phi$ determines a $\Spin(7)$ structure on $\Oc$.\footnote{In fact, $\Phi$ determines the Octonionic structure, a metric and an orientation on $\R^8$.}  The structure equations for the $\Spin(7)$-principal bundle \\
\centerline{\xymatrix{\Spin(7) \ar[r] &  {\Oc} \times \Spin(7) \ar[d]^{\bx} \\ & {\Oc}  }}\\  
 take the standard form
\be\label{eq:spin7streq}
\begin{split}
\ed \bx &= e_i \om_i\\
\ed e_i &= e_j \om_{ji}\\
\ed \om_{i}&= - \om_{ij} \w \om_{j}\\
\ed \om_{ij}&=-\om_{ik} \w \om_{kj}
\end{split}
\ee
where the $\om_{ij}$ satisfy the $\spin(7)$-relations given in equation \eqref{eq:spin7rel} of Appendix \ref{sec:spin7}. 

\begin{defn}
A smooth immersion $f:M^4 \to \Oc$ is a Cayley $4$-fold if 
\be\label{eq:cayley}
f^*(\psi_m)=0.
\ee
\end{defn} 
Cayley $4$-folds are calibrated by $\Phi$ and this implies that they are oriented.  Neither the perspective of calibrations nor the orientability of Cayley $4$-folds is needed in this article and so we refer the reader to \cite{hl} for the original definition of Cayley $4$-folds in terms of calibrations.  This article is concerned with the r-oriented Cayley cones ruled by $2$-planes, which are now defined.  

\begin{defn}
$M^4 \subset \Oc$ admits a \emph{smooth ruling by $2$-planes} (or is \emph{$2$-ruled}) if there exists a smooth surface $\Sigma$ and a smooth map $\pi:M\rightarrow \Sigma$ such that for all $\sigma \in \Sigma$, $\pi^{-1}(\sigma)$ is a $2$-plane in $\Oc$.  Such a triple $(M,\pi,\Sigma)$ is said to be a \emph{$2$-ruled} 4-fold.  If there exists a continuous choice of orientation for $\pi^{-1}({\sigma})$ then $M$ is said to be \emph{r-oriented}.
\end{defn}

Using parallel translation, each $2$-plane in $\{ \pi^{-1}(\sigma):\sigma \in \Sigma \}$ can be translated so that it contains the origin.  Thus $\pi^{-1}$ defines a map 
\be
\gamma:\Sigma \to \gro.
\ee
where $\gro$ is the Grassmannian of oriented two planes in $\Oc$.  The triple $(M,\pi,\Sigma)$ will be referred to as \emph{nondegenerate} if the associated map $\gamma$ is an immersion.  The $2$-dimensional family of $2$-planes $\gamma(\Sigma)$ forms a cone $M_0 \subset \Oc$.
\begin{defn}
$M_0$ is the \emph{asymptotic cone} of the {$2$-ruled} 4-fold $(M,\pi,\Sigma)$.  
\end{defn}
For a justification of the term \emph{asymptotic cone}, see Definition 3.2 in \cite{lotayca} and the surrounding text.  In Section \ref{sec:deformations} it will also be pointed out how each $2$-ruled Cayley $4$-fold comes in a one parameter family that has its asymptotic cone as a singular limit.
 
We now restrict our attention to the case in which $(M,\pi,\Sigma)$ is already a cone.  In Section \ref{sec:deformations} we return to the general case in order to prove that nondegenerate $2$-ruled Cayley cones are the asymptotic cones of smooth $2$-ruled Cayley $4$-folds whenever their associated surface $\gamma:\Sigma \to \gro$ satisfies a topological condition.

For a $2$-ruled cone $(M,\pi,\Sigma)$ the Cayley equations \eqref{eq:cayley} reduce to equations on the surface $\Sigma$.    A $2$-ruled cone in $\Oc$ can locally be parametrized on an open set $U \subset \Sigma$ as
\be
\bx: \R^2 \times U \to \Oc
\ee
\be
\bx (r_1,r_2,\sigma)=r_1 e_1 +r_2 e_2
\ee
where $e_1$ and $e_2$ are the first two legs of a $\Spin(7)$-adapted frame $e_i$.  Such a parametrization is always possible because $\Spin(7)$ acts transitively on oriented orthonormal pairs of vectors in $\Oc$.  This $4$-fold will be Cayley if $\bx^*(\psi_m)=0$ for $m=1 \ldots 7$.   One computes that
\be
\bx^*(\psi_m)=12 \,\ed r_1 \w \ed r_2 \w \psi_m(e_1,\,e_2,\,r_1e_i\omega_{i1}+r_2e_j\omega_{j2},\,r_1e_k\omega_{k1}+r_2e_l\omega_{l2}).
\ee
Expanding in $r_1$ and $r_2$ leads to the following system on $\Sigma$: 
\be\label{eq:cayleyeq}
\begin{split}
&\psi_m(e_1,\,e_2,\,e_i,\,e_k)\omega_{i1}\w\omega_{k1}=0\\
&\psi_m(e_1,\,e_2,\,e_i,\,e_l)\omega_{i1}\w\omega_{l2}+ \psi_m(e_1,\,e_2,\,e_j,\,e_k)\omega_{j2}\w\omega_{k1}=0\\
&\psi_m(e_1,\,e_2,\,e_j,\,e_l)\omega_{j2}\w\omega_{l2}=0
\end{split}
\ee 
For any such adapted coframe, define the complex forms
\begin{align}
\zeta_{3}&=\omega_{31}+i\omega_{41}\;\;\;\;\zeta_{4}=\omega_{32}+i\omega_{42} \nonumber \\ 
\zeta_{6}&=\omega_{61}-i\omega_{71}\;\;\;\;\zeta_{7}=\omega_{62}-i\omega_{72} \label{eqn:zetas}\\
\zeta_{5}&=\omega_{51}-i\omega_{81}\;\;\;\;\zeta_{8}=\omega_{52}-i\omega_{82}.\nonumber
\end{align}
All of the $2$-forms that appear in Equation \eqref{eq:cayleyeq} can be written as linear combinations of the real and imaginary parts of $\zeta_i \w \zeta_j$.  This indicates that, once the auxiliary directions of a $2$-ruled Cayley cone are boiled off, one is left with a geometry of pseudoholomorphic curves, which we now explain.   

\begin{defn} 
Let $(X^{2n},J)$ be a smooth manifold with a smooth almost complex structure and $\Sigma$ a smooth surface.  A \emph{pseudoholomorphic curve} is a smooth map
\be
\gamma:\Sigma \to (X,J)
\ee
for which $\gamma_*(T_{\sigma}\Sigma)$ is a complex line in $(T_xX,J_x) \cong \C^n$ whenever it is of real dimension two. The image $\gamma(\Sigma)$ inherits a complex structure and so $\Sigma$ acquires the structure of a Riemann surface. 
\end{defn}

The linear map $J$ (which also acts on $T^*_x X$) has eigenvalues $\pm i$ and so the module of smooth $\C$-valued $1$-forms $\Omega^1_{\C}( X)$ splits into the $+i$ and $-i$ eigenbundles $\Omega^{(1,0)}(X)$ and $\Omega^{(0,1)}(X)$, each with fiber isomorphic to $\C^n$.  Alternatively, one can define an almost complex structure on $X$ be choosing a splitting
\be\label{eqn:splitting}
\Omega^1_{\C}(X)=\Omega^{(1,0)}(X)\oplus \Omega^{(0,1)}(X)
\ee
such that $\ol{\Omega^{(1,0)}(X)} = \Omega^{(0,1)}(X)$.   An immersion of a (real) surface $\gamma:\Sigma \to X$ will be a pseudoholomorphic curve if and only if $\gamma^*(\alpha \w \beta)=0$ for all $\alpha,\beta \in \Omega^{(1,0)}(X)$, since this is equivalent to the tangent space of $\gamma(\Sigma)$ spanning a complex line in $T^{(1,0)}_{\gamma(\sigma)}X$.  It is this characterization of a pseudoholomorphic curve that we will use most often.

Similarly, for any real vector bundle $V$ with a complex structure, one can define the $(1,0)$-bundle inside $V \otimes \C$.  It will be denoted as $V'$ and $V''$ will denote the $(0,1)$-bundle.

In the case at hand, define a form $\alpha \in \Omega^1_{\C}(\gro)$ to be of type $(1,0)$ if for some (and thus any) $\Spin(7)$-adapted coframe  $\om_{ij}$, the $1$-form $\alpha$ is a section of the subbundle spanned by the $\zeta_i$ defined in Equation \eqref{eqn:zetas}.   This almost complex structure is $\Spin(7)$-invariant:

\newpage
\begin{prop}
\mbox \\
\begin{itemize}
\item{$\gro=\Spin(7)/\Un(3)$ and the $\Un(3)$ structure defines a $\Spin(7)$-invariant, nonintegrable almost complex structure}
\item{nondegenerate $r$-oriented $2$-ruled Cayley cones are equivalent to immersed pseudoholomorphic curves in $\gro$}
\end{itemize}
\end{prop}
\begin{proof}
The homogeneous structure is derived in \cite{hl}.  The nonintegrability of the $\Spin(7)$-invariant almost complex structure follows from the existence of the $\ol{\theta}_i \w \ol{\theta}_j$ terms in the structure equations derived below, Equation \eqref{eq:streq}.  This almost complex structure is referred to as $J_2$ in \cite{salamon} and the nonintegrability is given in Proposition 3.2 of \cite{salamon}. 

The rest of the proof is similar to that of the analogous result for $2$-ruled coassociative cones, Proposition 7.2 of \cite{fox}, so here we only sketch the proof of the second statement.  The Cayley conditions reduce to the vanishing of certain $2$-forms on the surface $\gamma(\Sigma) \subset \gro$ (Equation \eqref{eq:cayleyeq}) and these $2$-forms are the real and imaginary parts of $(2,0)$-forms for the almost complex structure on $\gro$ that is defined by declaring the $\zeta_i$ of Equation \eqref{eqn:zetas} to be of type $(1,0)$.  This shows that pseudoholomorphic curves in $\gro$ define $r$-oriented $2$-ruled Cayley cones.  The other direction is a little bit more involved.  Suppose we have an $r$-oriented $2$-ruled Cayley cone.  If there exists a coframe so that none of the $\zeta_i$ vanish on the surface $\Sigma$, then the Cayley equations \eqref{eq:cayleyeq} imply that $\zeta_i \w \zeta_j$ vanish on the surface, and thus it is a pseudoholomorphic curve on $\gro$. If some $\zeta_i$ vanishes on $\Sigma$, then the vanishing  of $\ed \zeta_i$ must be used to show that $\zeta_j \w \zeta_k=0$ on $\Sigma$ for all $j$ and $k$.  
\end{proof}
\noindent The nonintegrable $\Spin(7)$-invariant almost complex structure is closely related to an integrable one, as will be described in Section \ref{sec:fibration}.    

There is another perspective on pseudoholomorphic curves in $\gro$ which we will point out but not elaborate upon.  The almost complex structure used on $\gro=\Spin(7)/\Un(3)$ arises from its $3$-symmetric space structure.  Any $k$-symmetric space acts as the target for primitive maps \cite{bp} and in the case at hand the primitive maps are the same as the pseudoholomorphic curves.  Primitive maps are a natural generalization of harmonic maps and have a similarly elaborate theory in case that the domain is $\C$, $\s{2}$, or $T^2$.  This article is aimed at exploiting the exceptionally rich geometric structure of $\gro$ in order to obtain results about primitive maps from higher-genus domains. 
\end{section}

\begin{section}{The tautological bundles $\mch$ and $\mce$}\label{sec:taut1}
On $\gro$ there is a tautological complex line bundle $\mce$ and a tautological complex $3$-plane bundle $\mch$\\ 
\centerline{\xymatrix{ {\C \oplus \C^3} \ar[r] & {\mce} \oplus \mch \ar[d]^{(\pi_{\mce} ,\pi_{\mch})}\\ & \gro}}\\
which are defined as follows.  For $E \in \gro$ let $H=E^{\perp}$.  Then
\be
\pi_{\mce}^{-1}(E)=E;\;\;{\text and}\;\;\pi_{\mch}^{-1}(E)=H
\ee

The transitive action of $\Spin(7)$ on $\gro$ allows the frame to be chosen so that $E=e_1 \w e_2$.  Let
\be
f_0=e_1-ie_2,
\ee
and $f=(f_1,f_2,f_3)$ where 
\begin{align*}
f_1&=e_3-ie_4\\
f_2&=e_6+ie_7\\
f_3&=e_5+ie_8.
\end{align*}
Declaring the $f_0,f_i$ to be of type $(1,0)$ defines a complex structure on $\mce \otimes \C$ and $\mch \otimes \C$.    

We now calculate the structure equations for the tautological bundles $\mce$ and $\mch$.  A unitary framing of $\mce$ (of length $2$) is given by $f_0$ and $f=(f_1,f_2,f_3)$ is a unitary framing  of $\mch$ of length $2$.  The $\spin(7)$-structure equations \eqref{eq:spin7streq} can be organized as 
\begin{align}\label{eq:dg}
\ed (f_0,f) = (f_0, f)\bp -\tr(\kappa)& -\frac{1}{2}\trp \bar{\theta}_{ev} \\ \frac{1}{2}\theta_{ev} & \kappa \ep + (\fb_0,\fb) \bp 0& -\frac{1}{2}\trp\bar{\theta}_{od} \\ \frac{1}{2}\bar{\theta}_{od} & -\frac{i}{2} [\theta_{od}]\ep 
\end{align}
where the symbols have the following definitions:
\begin{align*}
\theta_{ev}&=\trp \bp \theta_2 & \theta_4 & \theta_6 \ep\\
\theta_{od}&=\trp \bp \theta_1 & \theta_3 & \theta_5 \ep
\end{align*}
\begin{align*}
\theta_{1}&=\zeta_{3}+i\zeta_{4}\;\;\;\;\theta_{2}=\zeta_{3}-i\zeta_{4}\\
\theta_{3}&=\zeta_{6}+i\zeta_{7}\;\;\;\;\theta_{4}=\zeta_{6}-i\zeta_{7}\\
\theta_{5}&=\zeta_{5}+i\zeta_{8}\;\;\;\;\theta_{6}=\zeta_{5}-i\zeta_{8}
\end{align*}
\begin{align*}
\kappa&=\bp i\om_{43} & -(\om_{63}+i\om_{64})-\frac{i}{2}\bar{\theta}_5 & -(\om_{53}+i\om_{54})+\frac{i}{2}\bar{\theta}_3\\
(\om_{63}-i\om_{64})-\frac{i}{2}{\theta}_5&-i\om_{76}& (\om_{65}-i\om_{75})-\frac{i}{2}\bar{\theta}_1\\
(\om_{53}-i\om_{54})+\frac{i}{2}{\theta}_3 & -(\om_{65}+i\om_{75})-\frac{i}{2}{\theta}_1&i(\om_{76}-\om_{43}-\om_{21}) \ep
\end{align*}
\noindent and satisfies 
\begin{align*}
{}^*\kappa&=-\kappa,
\end{align*}
and for any vector $v= \trp \bp v_1 & v_2 & v_3\ep$
\be
[v]=\bp 0&v_3&-v_2 \\ -v_3 &0 & v_1 \\ v_2 & -v_1 &0\ep.
\ee
The map $v \to [v]$ gives an explicit isomorphism $\C^3 \to \Lambda^2 \C^3$.

The form of the structure equations \eqref{eq:dg} modulo ${\theta}_{od},{\theta}_{ev},\ol{\theta}_{od},\ol{\theta}_{ev}$ indicate that the stabilizer of $E$ is isomorphic to $\Un(3)$ and that it acts on $H$ in the standard way and on $E$ as $\det^{-1}$. 

Define the metric compatible connections $\nabla^{\mch}$ and $\nabla^{\mce}$ as
\begin{align*}
\nabla^{\mch}f&=\kappa f\\
\nabla^{\mce} f_0&=-{\tr(\kappa)} f_0
\end{align*}
Both $\mch'$ and $\mce'$ pull back to any pseudoholomorphic curve to be Hermitian using these connections.  The curvature equations for these connections will also be used.  For $\mce$:
\be
\ed (-\tr(\kappa))=\frac{1}{4}(\theta_{2} \w \bar{\theta}_{2} +\theta_{4} \w \bar{\theta}_{4} +\theta_{6} \w \bar{\theta}_{6} -\theta_{1} \w \bar{\theta}_{1} -\theta_{3} \w \bar{\theta}_{3} -\theta_{5} \w \bar{\theta}_{5} )
\ee
and for $\mch$:
\be
\ed \kappa+ \kappa \w \kappa =\frac{1}{4}\Omega, 
\ee
where $\Omega$ is skew-hermitian with components
\begin{center}
\begin{tabular}{ll}
$\Omega_{11} =\;\;\; \theta_{1} \w \bar{\theta}_{1}- \theta_{3} \w \bar{\theta}_{3}- \theta_{5} \w \bar{\theta}_{5}+ \theta_{2} \w \bar{\theta}_{2}$& $\;\;\;\;\;\Omega_{21} = 2\theta_{3} \w \bar{\theta}_{1} +\theta_{4} \w \bar{\theta}_{2}$\\
$\Omega_{22} = -\theta_{1} \w \bar{\theta}_{1}+ \theta_{3} \w \bar{\theta}_{3}- \theta_{5} \w \bar{\theta}_{5}+ \theta_{4} \w \bar{\theta}_{4}$& $\;\;\;\;\;\Omega_{31} = 2\theta_{5} \w \bar{\theta}_{1} +\theta_{6} \w \bar{\theta}_{2}$\\
$\Omega_{33} = -\theta_{1} \w \bar{\theta}_{1}- \theta_{3} \w \bar{\theta}_{3}+ \theta_{5} \w \bar{\theta}_{5}+ \theta_{6} \w \bar{\theta}_{6}$ &$\;\;\;\;\; \Omega_{32} = 2\theta_{5} \w \bar{\theta}_{3} +\theta_{6} \w \bar{\theta}_{4}$.\\
\end{tabular}
\end{center}
\end{section}

\begin{section}{Structure equations for the tangent bundle of $\gro$}\label{sec:tangstr}
We next derive the structure equations for the $\Un(3)$-structure on $\gro$.  Let
\be
\bfp= 
e_1 \w e_2=-\frac{i}{2}f_0 \w \bar{f_0}:\Spin(7) \to \gro.
\ee
We differentiate:
\begin{align*}
\ed \bfp &=\frac{1}{4}\left[ if_1 \w f_0 \theta_1 - i \bar{f}_{1} \w \bar{f}_0 \bar{\theta}_1\right]
-\frac{1}{4}\left[ if_1 \w \bar{f}_0 \theta_2 - i \bar{f}_{1} \w {f}_0 \bar{\theta}_2 \right]\\
&+\frac{1}{4}\left[ if_{2} \w f_0 \theta_3 - i \bar{f}_{2} \w \bar{f}_0 \bar{\theta}_3\right]
-\frac{1}{4}\left[ if_{2} \w \bar{f}_0 \theta_4 - i \bar{f}_{2} \w {f}_0 \bar{\theta}_4 \right]\\
&+\frac{1}{4}\left[ if_{3} \w f_0 \theta_5 - i \bar{f}_{3} \w \bar{f}_0 \bar{\theta}_5\right]
-\frac{1}{4}\left[ if_{3} \w \bar{f}_0 \theta_6 - i \bar{f}_{3} \w {f}_0 \bar{\theta}_6 \right].
\end{align*}

Thus the $(1,0)$-vectors are spanned by $f_i \w f_0$ and $f_i \w \bar{f}_0$ for the dual $(1,0)$-coframe $\theta_i$.  This shows that the $(1,0)$-subbundle of the complexified tangent bundle $T\gro \otimes \C$ splits,
\be
T^{(1,0)}\gro = \vooz \oplus \vtoz,
\ee
where each bundle $\mcv'_i$ has fiber $\C^3$.  By comparing the respective bases of the $(1,0)$-vectors in $\mch$, $\mce$, $\mcv_1$ and $\mcv_2$ we find that 
\begin{align*}
\mcv_1 &\cong \mch \otimes \mce\\
\mcv_2 &\cong \mch \otimes \mce^{*}.
\end{align*}

The structure equations for a coframe of $\gro$ are 
\be\label{eq:streq}
\ed \begin{bmatrix} \theta_{od}  \\ \theta_{ev} \end{bmatrix}  = -
\begin{bmatrix}
\Psi & 0 \\ 0 &  \tilde{\Psi}
\end{bmatrix}
\w \begin{bmatrix} \theta_{od} \\ \theta_{ev} \end{bmatrix} + \begin{bmatrix}\tau_{od}\\ \tau_{ev} \end{bmatrix}
\ee
where 
\begin{align*}
\Psi&=\kappa -\tr(\kappa)\\
\tilde{\Psi}&=\kappa +\tr(\kappa)
\end{align*}
and
\begin{align*}
\tau_{1} &= \frac{i}{2} \left( \bar{\theta}_{4} \w \bar{\theta}_{5}+ \bar{\theta}_{3} \w \bar{\theta}_{6} \right) \hspace{2cm} \tau_{2} = {i} \bar{\theta}_{3} \w \bar{\theta}_{5}    \\
\tau_{3} &= \frac{i}{2} \left( \bar{\theta}_{6} \w \bar{\theta}_{1}+ \bar{\theta}_{5} \w \bar{\theta}_{2} \right) \hspace{2cm} \tau_{4} = {i} \bar{\theta}_{5} \w \bar{\theta}_{1} \\ 
\tau_{5} &= \frac{i}{2} \left( \bar{\theta}_{2} \w \bar{\theta}_{3}+ \bar{\theta}_{1} \w \bar{\theta}_{4} \right) \hspace{2cm} \tau_{6} = {i} \bar{\theta}_{1} \w \bar{\theta}_{3}.  
\end{align*}
\end{section}

\begin{section}{The fibration $\gro \to \s{6}$ and the distributions $\mcv_1$ and $\mcv_2$}\label{sec:fibration}
The results of this section are well known in the twistor theory literature.  For example, see \cite{eellssalamon,salamon,rawnsley,burstallrawnsley} and chapter 9 of \cite{lawsonmichelson}.

The splitting $T\gro=\mcv_1 \oplus \mcv_2$ is induced by the fibration\\
\centerline{\xymatrix{\cp{3} \ar[r] & \gro \ar[d]^{\mcj} \\ & \s{6}} \label{eq:ssfib}}\\
\\
and plays a fundamental role in the geometry of pseudoholomorphic curves in $\gro$.  In this section the fibration and its most basic implications will be described.  In Section \ref{sec:sixsphere} a more detailed investigation of the fibration will locate the truly \emph{holomorphic} aspects of the \emph{pseudoholomorphic} curves in $\gro$.

Using the exceptional isomorphism 
\be\label{eq:except}
\Spin(6) \cong \SU(4),
\ee
the map $\mcj$ can be defined as follows.  Each $2$-plane $E \in \gro$ has stabilizer  $\Un(3) \subset \Spin(7)$ and thus is given by a coset $g \Un(3) \in \Spin(7)/\Un(3)$.  There is a unique $\SU(4) \subset \Spin(7)$ containing that stabilizer.   The map is 
\begin{align*}
\mcj&:\Spin(7)/\Un(3) \to \Spin(7)/\SU(4)\\
\mcj&:g\Un(3) \mapsto g\SU(4).
\end{align*}
The exceptional isomorphism (Equation \eqref{eq:except}) provides the identification 
$$\Spin(7)/SU(4)=\Spin(7)/\Spin(6)=\s{6}.$$
This gives $\s{6}$ the interpretation as the space of possible reductions from the $\Spin(7)$-structure to an $\SU(4)$-structure. 

It will be useful to understand how the fibration is manifested in the infinitesimal geometry.  
\begin{lem}
\mbox \\
\begin{itemize}
\item{The distribution $\mcv_1$ is not integrable.  The map   $\mcj_*:\mcv_1 \to T\s{6}$ is an isomorphism.}
\item{The distribution $\mcv_2$ is integrable and its maximal leaves are the $\cp{3}$-fibers of $\mcj$.}
\end{itemize}
\end{lem}

\begin{proof}
The first distribution is defined by the ideal 
\be
\mci_1=\langle \Re(\theta_{ev}), \Im(\theta_{ev})\rangle
\ee
and the second by the ideal 
\be
\mci_2=\langle \Re(\theta_{od}), \Im(\theta_{od})\rangle.
\ee
Structure equations \eqref{eq:streq} show that the first ideal is not Frobenius while the second is.

The maximal leaves of $\mci_2$ are six dimensional manifolds.  On a leaf the structure equations from \eqref{eq:dg} become
\be
\ed (f_0,f) = (f_0, f)\bp -\tr(\kappa)& -\frac{1}{2}\trp \bar{\theta}_{ev} \\ \frac{1}{2}\theta_{ev} & \kappa \ep. 
\ee 
Thus they are $\SU(4)$-orbits with stabilizer $\Un(3)$, i.e., they are isomorphic to $\cp{3}$.  Since there is a unique $\SU(4) \subset \Spin(7)$ containing any given $\Un(3)$, these leaves must be the fibers of the map $\mcj$.

One checks that $\ker(\mcj_*)=\mcv_2$ and then the isomorphism $\mcj_*:\mcv_1 \to T\s{6}$ follows from a dimension count.
\end{proof}

The fibration distinguishes two special types of pseudoholomorphic curves in $\gro$: those that are in a single fiber of $\mcj$ (tangent to $\mcv_2$) and those that are \emph{horizontal} (tangent to $\mcv_1$).  It is clear that the first type consists of actual algebraic curves in $\cp{3} \subset \gro$.  Let $P \in \s{6}$ be the image under $\mcj$ of such a pseudoholomorphic curve.  It defines an $\SU(4)$-structure on $\Oc$ and thus a complex structure $J$.  The pseudoholomorphic curves that are contained in a fiber $\mcj^{-1}(P)$ correspond to the $2$-ruled Cayley cones in $\Oc$ that are actually holomorphic cones in $\C^4 \cong \Oc_J$.  

We will now see that the second type also consists of algebraic curves.  Let  
\be
Q_6=\{[z] \in \cp{7} \;|\; z_1^2+ \ldots + z_8^2=0 \} \subset \cp{7}
\ee
denote the six quadric, which is diffeomorphic to $\gro$ via the map 
\begin{align*}
&\gro \to Q_6\\
&e_1 \w e_2 \mapsto [e_1 - ie_2]. 
\end{align*}
The $\Spin(7)$-invariant almost complex structure is not integrable and so does not induce the holomorphic structure on $Q_6$.  However the two almost complex structures do agree on $\mcv_1$.  
\begin{lem}\label{lem:dist}
\mbox{}
\begin{itemize}
\item{The complex $3$-plane distribution $\mcv'_1$ is holomorphic with respect to the integrable complex structure on $Q_6$.}
\item{The complex $3$-plane distribution $\ol{\mcv'_2}$ is holomorphic with respect to the integrable complex structure on $Q_6$.}
\end{itemize}
\end{lem}

To prove Lemma \ref{lem:dist} we need to relate the holomorphic structure on $Q_6$ and the $\Spin(7)$-invariant almost complex structure.  In an $\SO(8)$-adapted coframe (and so also in a $\Spin(7)$-adapted coframe) the $(1,0)$-forms of the holomorphic structure on $Q_6$ are spanned by
\be\label{eqn:xis}
\xi_j=\om_{j1}+i\om_{j2}
\ee
and satisfy 
\be
\ed \xi_j = -i\om_{21} \w \xi_j -\om_{jk} \w \xi_k \;\;\;(j,k >2).
\ee
The closure of the algebraic ideal $\langle \xi_i \rangle$ under exterior differentiation implies the integrability of the corresponding complex structure.  Equations \eqref{eqn:xis} and \eqref{eqn:zetas} allow one to relate the $\xi_i$ and the $\zeta_i$. 

\begin{proof}
We begin with $\mcv'_1$.  The vanishing of the $\theta_{2j}$ is equivalent to
\begin{align*}
\xi_{4}&=-i\xi_{3}\\
\xi_{7}&=i\xi_{6}\\
\xi_{8}&=i\xi_{5}
\end{align*}
which define a holomorphic distribution on $Q_6$.  One checks that the complex structures on $\mcv_1'$ and $Q_6$ agree.  

Now consider $\ol{\mcv'_2}$. The vanishing of the $\theta_{2j-1}$ is equivalent to 
\be\label{eq:holdist2}
\begin{split}
\xi_{4}&=i\xi_{3}\\
\xi_{7}&=-i\xi_{6}\\
\xi_{8}&=-i\xi_{5}
\end{split}
\ee
which define a holomorphic distribution on $Q_6$.  One can check that the complex structure of $\ol{\mcv'_2}$ agrees with that of $Q_6$.  
\end{proof}

\begin{cor}\label{cor:algebraiccurves}
Let $\gamma:\Sigma \to \gro$ be a pseudoholomorphic curve in $\gro$.
\begin{itemize}
\item{If it is tangent to $\mcv_1$, $\gamma(\Sigma)$ is algebraic in $Q_6$.}
\item{If it is tangent to $\mcv_2$ then $\overline{\gamma(\Sigma)}$ is algebraic in $\cp{3} \subset Q_6$.}
\end{itemize}
\end{cor}

\end{section}

\begin{section}{First fundamental forms}\label{sec:fff}
There are many invariants for pseudoholomorphic curves in $\gro$.  We now introduce two basic first order invariants.  Throughout this section we will make use of the natural identifications
\be
{\mathbb P}(\mcv'_i \otimes (T^*\Sigma)') \cong {\mathbb P}(\mcv'_i) \cong {\mathbb P}(\hoz ) 
\ee
that arise from the natural isomorphisms $\mcv_1 = \mch \otimes \mce$ and $\mcv_2 = \mch \otimes \mce^{-1}$.  Correspondingly, a (complex) line in any of these vector spaces will be freely identified with a line in the others.
   
Let 
\be
\gamma:\Sigma \to \gro
\ee
be a pseudoholomorphic curve.  Denote the tangent bundle $T\Sigma$ by $T$.  
\begin{defn}
The first fundamental forms of $\gamma$ are
\begin{align}
\ffo&=\gamma^*(f_0 \otimes f \otimes \thod) \in \Gamma(\gamma^*(\eoz \otimes \hoz ) \otimes T'^*)\\
\ffe&=\gamma^*(\fb_0 \otimes f \otimes \thev) \in \Gamma(\gamma^*(\ezo \otimes \hoz ) \otimes T'^*).
\end{align}
\end{defn} 

\begin{lem}
$\ffo$ and $\ffe$ are holomorphic sections and $\ffo$ \us{(} resp. $\ffe$\us{)} vanishes identically if and only if $\gamma(\Sigma)$ is tangent to $\mcv_2$ \us{(} resp. $\mcv_1$\us{)}. 
\end{lem}

The proof of the Lemma, which is just a calculation in local coordinates, takes the same shape as the proof of Lemma 4.2 in \cite{boc} or Lemma 3.1 in \cite{xu} and we direct the reader to either proof for details.  The holomorphicity of $\ffe$ is also a special case of Lemma 6.18 of \cite{burstallrawnsley}

As long as $(\Sigma,\gamma)$ is not tangent to one of the holomorphic distributions, $\ffo$ and $\ffe$ only vanish at isolated points.  Thus there exist holomorphic line subbundles of $\gamma^*(\mcv_{i}')$, of which the ${\rm I}_i$ are nonzero sections.  For later use (see Lemma \ref{lem:minsurf}, Section \ref{sec:highergenus}, and Section \ref{sec:deformations}) we define $R_1$ and $R_2$ to be the vanishing loci of $\ffo$ and $\ffe$.
\end{section}

\begin{section}{The relationship to minimal surfaces in $\s{6}$}\label{sec:sixsphere}
We now consider in more detail the implications of the fibration\\
\centerline{\xymatrix{\cp{3} \ar[r] & \gro \ar[d]^{\mcj} \\ & \s{6}}}\\
on the geometry of pseudoholomorphic curves in $\gro$.  The structure equations for $\s{6}$ are related to those of $\gro$ as follows.  Let
\be
\eta_i=\Re(\theta_i)\;\;\sigma_i=\Im(\theta_i)\;\;\p=\Re(\kappa)\;\;\psi=\Im(\kappa).
\ee
Then  
\be\label{eq:s6streq}
\ed \bp \eta_{od}\\ \sigma_{od} \ep =- \bp \phi+ \frac{1}{2}[\sigma_{ev}]& -\trp(\psi + \frac{1}{2}[\eta_{ev}]) \\ \psi + \frac{1}{2}[\eta_{ev}] & \phi- \frac{1}{2}[\sigma_{ev}] \ep \w \bp \eta_{od}\\ \sigma_{od} \ep 
\ee
are the structure equations for the standard $\SO(6)$-structure on $\s{6}$.

The fibration relates pseudoholomorphic curves in $\gro$ to minimal surfaces in $\s{6}$.
\begin{lem}\label{lem:minsurf}
Let $\gamma:\Sigma \to \gro$ be a pseudoholomorphic curve.  Then 
$$\phi=\mcj \circ \gamma:\Sigma \to \s{6}$$
is a minimal surface with a finite number of branch points at $R_1$, the zero locus of the first fundamental form of $\gamma$.   
\end{lem}
This is due to Salamon (\cite{salamon}, Theorem 3.5) so we only sketch a proof.

\begin{proof}
Locally choose an adapted coframe on the pseudoholomorphic curve such that 
\be\label{eq:3forms}
\theta_3=\theta_5=\theta_6=0.
\ee
On the one hand this implies that $\eta_3=\eta_5=\sigma_3=\sigma_5=0$ so that $\eta_1,\sigma_1$ are a coframe for $\phi(\Sigma) \subset \s{6}$.  In the standard way this allows one to express the second fundamental form in terms of certain connection forms from \eqref{eq:s6streq}.  On the other hand, differentiating the vanishing conditions in \eqref{eq:3forms} implies that $\theta_1$, $\theta_2$, $\theta_4$, $\kappa_{21}$, $\kappa_{31}$, $\kappa_{32}$ are all of type $(1,0)$.  This now implies that the second fundamental form of $\phi(\Sigma)$ is trace free.

The fact concerning branch points follows from the definition of $I_1$ and the fact that the real and imaginary parts of $\theta_{od}$ span the semi-basic forms for the fibration $\mcj:\gro \to \s{6}$.
\end{proof}

\begin{rem}
As we will see below, in general the map 
\be
\phi: \Sigma \to \s{6}
\ee
is a multiple-sheeted branched cover onto its image $\phi(\Sigma) \subset \s{6}$.  
\end{rem}

Lemma \ref{lem:minsurf} invites the question:  Is every minimal surface in $\s{6}$ the image of a pseudoholomorphic curve?  The minimal surfaces in $\s{6}$ locally depend on eight functions of one variable whereas the pseudoholomorphic curves in $\gro$ depend on ten functions of one variable, which suggests an affirmative answer.  Moreover this function count suggests that there should be two functions of one variable's worth of pseudoholomorphic curves sitting over any minimal surface, at least locally.  To understand this better it is useful to consider the fibration from a slightly different perspective.  

Let \\
\centerline{\xymatrix{ {\C}^4 \ar[r] & \spn \ar[d] \\ & \s{6} }}\\
be the spinor bundle on $\s{6}$ induced from the homogeneous structure 
$$\s{6}=\Spin(7)/\Spin(6).$$
The action of $\Spin(6)$ on the fibers is via the isomorphism $\Spin(6)=\SU(4)$ and the standard action of $\SU(4)$ on $\C^4$.

\begin{lem}
As fiber bundles over $\s{6}$, $\gro \cong {\mathbb P}(\spn)$.
\end{lem}
\begin{proof}
The preimage under $\mcj$ of a point $P \in \s{6}$ is an $\SU(4)$ orbit since $\stab(P)=\Spin(6) \cong \SU(4) \subset \Spin(7)$. Let $J$ be the complex structure on $\Oc$ defined by $P$.  The orbit is a $\cp{3}=\SU(4)/\Un(3)$ since the stabilizer of a point in $\gro$ is $\Un(3)$.  So any point in $\mcj^{-1}(P)$ is a complex line in $\C^4=\Oc_J$, which is the same as the fiber of ${\mathbb P}(\Delta)$.
\end{proof}

Now we begin to work out how the geometry of the spinor bundle captures the geometry of a minimal surface.  Let
\be
\phi:\Sigma \to \s{6}
\ee
be a branched minimal immersion of an \emph{oriented} surface.  This induces a conformal structure on $\Sigma$. Let $T$ denote the tangent bundle of $\phi(\Sigma)$, $T'$ the corresponding $(1,0)$ component of the complexified tangent bundle, and let $N$ denote the normal bundle of $\phi(\Sigma) \subset \s{6}$, a real rank-$4$ bundle.   

The product formula for Stiefel-Whitney classes implies that the normal bundle of an oriented surface in a six dimensional spin manifold has a spin structure.  Denote the principal spin bundle as $\Spin(N)=\Spin(N)^+ \times \Spin(N)^-$ and the associated spinor bundles as $\spn_{\pm}$.  A change in orientation on $\Sigma$ will interchange $\spn_+$ and $\spn_-$, so we can focus our attention on $\spn_+$ and similar results follow for $\spn_-$ by changing the orientation of $\Sigma.$
\begin{lem}
There is a decomposition
\be
\phi^*(\spn)=\mcw_+ \oplus \mcw_-
\ee
where 
\begin{align*}
\mcw_{+}&=\spn_{+} \otimes T'\\
\mcw_{-}&=\spn_{-} \otimes (T')^*
\end{align*}
and their principal bundles are respectively
\begin{align*}
\Spin^c(N)^{+}&=(\Spin(N)^{+} \times \Spin(T))/\Z_2\\
\Spin^c(N)^{-}&=(\Spin(N)^{-} \times \Spin(T))/\Z_2.
\end{align*}
\end{lem}
\begin{proof}
The decomposition 
\be
\phi^*(T\s{6})=T\Sigma \oplus N\Sigma
\ee
corresponds to the reduction of the structure group 
\begin{center}
\begin{tabular}{c}
$\Spin(6)$  \\
$\cup$  \\ 
$(\Spin(2) \times \Spin(3)^+ \times \Spin(3)^-)/\Z_2$. 
\end{tabular}
\end{center}
Using the exceptional isomorphisms this is equivalent to the reduction
\begin{center}
\begin{tabular}{c}
  $\SU(4)$ \\
 $ \cup$\\ 
 $S(\Un(2) \times \Un(2))$, 
\end{tabular}
\end{center}
which implies the splitting of the spinor bundle $\spn$ given above.
\end{proof}

\begin{lem}\label{lem:hollift}
Suppose that $\phi:\Sigma \to \s{6}$ is minimal and oriented.  Any holomorphic line subbundle $L \subset \mcw_{+}$ defines a \emph{pseudoholomorphic} lift
\be
\gamma_{(\phi,L)}:\Sigma \to \gro.
\ee
\end{lem}
\begin{proof}
The prescription defines a lift because $L \in \mathbb{P}(\mcw_+) \subset \mathbb{P}( \phi^*\spn)=\phi^*\gro$.  We must check that it is holomorphic and to do so we calculate locally.  In a neighborhood of any point there exists a coframe in which $\theta_3=\theta_5=0$ and then the minimality implies that $\theta_4$, $\theta_6$, $\kappa_{21}$, $\kappa_{31}$ are of type $(1,0)$ with respect to the holomorphic structure induced by $\theta_1$.  

Locally we have 
\be
\mcw_+=\C \cdot \{f_0,f_1 \} 
\ee
and this bundle has the induced connection
\begin{align*}
\nabla^+ (f_0,f_1)&=(f_0,f_1) \bp -\tr(\kappa)& -\frac{1}{2}\bar{\theta}_2 \\ \frac{1}{2}\theta_2 & \kappa_{11} \ep 
\end{align*}
A further frame adaptation allows $L=\C \cdot \{ f_0\}$.  It then follows that locally
\[
\ff_L=f_1 \otimes f^*_0 \otimes \theta_2,
\]
where $\ff_L$ is the second fundamental form of the line subbundle $L$.  In Appendix \ref{sec:sff} the second fundamental form of a line subbundle is defined and  Lemma \ref{lem:hlb} implies that $\theta_2$ is of type $(1,0)$ since $L$ is holomorphic. This implies that 
\be
\gamma_{(\phi,L)}=[f_0] \circ \phi:\Sigma \to \gro
\ee
is a pseudoholomorphic curve since $\theta_i \w \theta_j=0$ for all $i,j$.    
\end{proof}

\begin{rem}
Choosing a holomorphic line subbundle of $\mcw_+$ is locally equivalent to choosing a holomorphic map $\Sigma \to \cp{1}$ which depends on two functions of one variable.  This clarifies the function count made earlier. 
\end{rem}

\begin{cor}
The pseudoholomorphic curve $(\Sigma,\gamma_{(\phi,L)})$ projects to $\phi(\Sigma)$ in a one-to-one fashion and so every minimal surface in $\s{6}$ is the image of a pseudoholomorphic curve in $\gro$.  This is true even if $\phi$ has branch points.  
\end{cor}
\noindent This corollary is also a special case of a theorem of Rawnsley (\cite{rawnsley}, Theorem 9.10).

The results of this section are summarized in the following theorem, which decomposes every pseudoholomorphic curve in $\gro$ into a (possible branched) minimal surface in $\s{6}$ and holomorphic data.
\begin{thm}\label{thm:data}
Any pseudoholomorphic curve $\gamma:\Sigma \to \gro$ is equivalent to the following data: 
\begin{center}
\begin{tabular}{ll}
$\bullet \;\phi':\Sigma' \to \s{6} $ & a branched minimal immersion of an oriented surface\\
$\bullet \;\psi:\Sigma \to \Sigma'$ & a branched cover of Riemann surfaces\\
$\bullet \;\phi:\Sigma \to \s{6}$ & a branched minimal immersion of an oriented surface\\
$\bullet \;L$ & a holomorphic line subbundle of $ \phi^*(\mcw_+)$\\
$\bullet \;\gamma=\gamma_{(\phi,L)}$& where this map is defined as in Lemma \ref{lem:hollift}, \\
\end{tabular}
\end{center}
where $\phi=\phi' \circ \psi$ and $\phi(\Sigma')=\mcj \circ \gamma(\Sigma)$.
\end{thm}

\begin{rem}
The bulk of the geometry of a pseudoholomorphic curve in $\gro$ arises from a minimal surface in $\s{6}$.  Much is known about minimal spheres in $\s{6}$.  Minimal spheres are always isotropic and thus lift to $\gro$ to be tangent to $\mcv_1$.  They are therefore algebraic curves in $Q_6$ and it has been shown that the moduli space of minimal spheres of fixed degree $\deg$ is an algebraic variety of complex dimension $9+2\deg$.  The short survey by Bolton and Woodward \cite{bw} addresses these main points but is a little bit out of date by now.  

Less is known about minimal tori in $\s{6}$, though there is a spectral curve.  It would be nice to have an existence theorem for such tori along the lines of the work of Carberry and McIntosh \cite{c} \cite{cmci}.  

For higher-genus minimal surfaces it may be fruitful to mimic the twistor methods of Bryant \cite{boc}.  He shows that every compact Riemann surface admits a \emph{branched} immersion into $\s{6}=\G/\SU(3)$ as a null-torsion pseudoholomorphic curve.  The null-torsion condition is similar to the isotropic condition of \cite{bw}. It seems likely that Bryant's approach could be extended to prove the existence of branched minimal surfaces of every genus in $\s{6}$ that are not pseudoholomorphic for any $\G$-invariant almost complex structure.  These would all arise from pseudoholomorphic curves in $\gro$ that are tangent to $\mcv_1$. Xu \cite{xu} used this method to prove the existence of higher genus pseudoholomorphic curves in the nearly K\"ahler $\cp{3}$.         
\end{rem}

Without going into detail I want to describe an interesting degenerate case.  From this degenerate case one recovers the surprising duality between minimal surfaces in $\s{5}$ and ruled minimal Lagrangian $3$-folds in $\cp{3}$, which appeared in the work of Bolton and Vrancken \cite{bv}.  They work directly with a moving frame instead of using the $2$-ruled special Lagrangian cones as a link between the two geometries.   At the most basic level, this equivalence is a consequence of the exceptional isomorphism $\Spin(6) \cong \SU(4)$.

Choosing a point $P \in \s{6}$ defines both the perpendicular five sphere $\iota:\s{5}_P \to \s{6}$ (if $P$ is a pole, $\s{5}_P$ is the equator) and an $\SU(4)$-structure on $\Oc$.  Let $(\Omega,\Upsilon)$ be the K\"ahler and holomorphic volume forms for the $\SU(4)$-structure.  One finds that $\Phi=\frac{\Omega^2}{2}+\Re(\Upsilon)$.  If $\Re(\Upsilon)$ vanishes on a $2$-ruled Cayley cone then it is a $2$-dimensional holomorphic cone in $\C^4$ and this is equivalent to the corresponding pseudoholomorphic curve being contained in the fiber $\cp{3}_P=\mcj^{-1}(P)$.  

If $\frac{\Omega^2}{2}$ vanishes on a $2$-ruled Cayley cone then it is a special Lagrangian cone and the associated pseudoholomorphic curve is contained in $\mcj^{-1}(\s{5}_P)$.  In fact $\mcj^{-1}(\s{5}_P)=\gris$, where $\gris$ is the Grassmannian of oriented isotropic real $2$-planes in $\C^4$.  The image of the pseudoholomorphic curve in $\s{5}_P$ is necessarily minimal.  One can check that there is at most a discrete family of pseudoholomorphic curves in $\gris$ that project to the same minimal surface in $\s{5}_P$ even locally.  

Given a minimal surface in $\s{5}_P$ one can lift it to a pseudoholomorphic curve in $\gris \subset \gro$ and this defines a $2$-ruled special Lagrangian cone $L \subset \C^4$.  Then the link $M^3=L^4 \cap \s{7}$ projects down to $\cp{3}=\s{7}/\s{1}$ to be a minimal Lagrangian submanifold that is fibered by geodesic $\s{1}$'s.  This process can clearly be reversed and this gives a correspondence between minimal surfaces in $\s{5}$ and minimal Lagrangian $3$-folds (ruled by geodesic $\s{1}$'s) in $\cp{3}$.          
\end{section}

\begin{section}{The degree}
We quickly review the definition of the degree of a pseudoholomorphic curve in $\gro$.  It will be used in Sections \ref{sec:highergenus} and \ref{sec:deformations}.

The $(1,1)$-form on $\gro$ is 
\be
\omega=\frac{i}{2} [\theta_{1} \w \bar{\theta}_{1}+\theta_{3} \w \bar{\theta}_{3} +\theta_{5} \w \bar{\theta}_{5} +\theta_{2} \w \bar{\theta}_{2}  +\theta_{4} \w \bar{\theta}_{4} +\theta_{6} \w \bar{\theta}_{6} ].
\ee
The splitting $\Omega^{(1,0)}(\gro) \cong {\mcv_1^*}' \oplus {\mcv_2^*}'$ provides a decomposition 
\be
\omega=\om_1 +\om_2
\ee
into the odd and even halves.  These $(1,1)$-forms are not closed:
\be
\ed \om_1 = \ed \om_2 = -\Re(\theta_{2} \w \theta_{3} \w \theta_{5} +\theta_{1} \w \theta_{4} \w \theta_{5} +\theta_{1} \w \theta_{3} \w \theta_{6} ),
\ee
but their difference is.  One can check that $\om_1-\om_2$ is the standard K\"ahler form on $\gro$ induced by the diffeomorphism $\gro \cong Q_6$.   

A topological measure of the complexity of a pseudoholomorphic curve in $\gro$ is its \emph{degree}.  
\begin{defn}
The degree of a pseudoholomorphic curve $\gamma:\Sigma \to \gro$ is defined to be
\be
\ed_{\gamma}:=\int_{\Sigma} \gamma^*(c_1(\mce))
\ee
where $c_1(\mce)$ is the first Chern class of $\mce$.
\end{defn}
The connection on $\mce$ is $-\tr(\kappa)$, so that its first chern class is 
\begin{align}
c_1(\mce)&=[\frac{i}{2\pi} \ed (-\tr(\kappa))]\\
&=\frac{1}{4\pi}\frac{i}{2} [(\theta_{1} \w \bar{\theta}_{1}+\theta_{3} \w \bar{\theta}_{3} +\theta_{5} \w \bar{\theta}_{5}) -(\theta_{2} \w \bar{\theta}_{2}  +\theta_{4} \w \bar{\theta}_{4} +\theta_{6} \w \bar{\theta}_{6}) ]\\
\end{align}
or
\be
c_1(\mce)=\frac{1}{4\pi}(\om_1-\om_2).
\ee
\end{section}

\begin{section}{The existence of higher genus immersed curves}\label{sec:highergenus}
The characterization of pseudoholomorphic curves in $\gro$ provided by Theorem \ref{thm:data} leads to a simple method for producing immersed higher genus pseudoholomorphic curves from those of lower genus, or simply from a minimal surface in $\s{6}$.  Begin with a (possibly branched) minimal surface 
\be
\phi':\Sigma' \to \s{6}.
\ee
Now choose a higher-genus Riemann surface $\Sigma$ that admits a branched cover $\psi:\Sigma \to \Sigma'$ and let $\phi=\phi' \circ \psi$.  This defines a branched minimal immersion 
\be
\phi:\Sigma \to \s{6}
\ee
that multiply covers itself.  This induces the rank two Hermitian vector bundle $\mcw_+$ over $\Sigma$.  To lift $\Sigma$ to a pseudoholomorphic curve in $\gro$ requires choosing a holomorphic line subbundle of $\mcw_+$. Such subbundles are bountiful.   
\begin{thm}\label{thm:immersed}
Let $(\Sigma',\phi')$ be a minimal immersion, $\psi:\Sigma \to \Sigma'$ a branched cover of Riemann surfaces, and $\phi=\phi' \circ \psi:\Sigma \to \s{6}$.
\begin{itemize}
\item{There exist (possibly branched) immersions $(\Sigma,\gamma)$ as pseudoholomorphic curves such that $\phi=\mcj \circ \gamma$ and for which $(\Sigma,\gamma)$ has arbitrarily large negative degree.}
\item{There exist \emph{immersions} $\gamma:\Sigma \to \gro$ as pseudoholomorphic curves such that $\phi=\mcj \circ \gamma$.}
\end{itemize}
\end{thm}
\begin{rem}
One could choose $\Sigma'=\cp{1}$ and then every compact Riemann surface becomes a branched cover of it by choosing a meromorphic function on the Riemann surface.  Therefore every compact Riemann surface appears as an \emph{immersed} pseudoholomorphic curve in $\gro$.
\end{rem}

The proof of Theorem \ref{thm:immersed} follows from the following sequence of lemmas.  From the discussion prior to the statement of the theorem, we must find a holomorphic line subbundle of $\mcw_+$ with the desired properties.  We will first  point out why every Hermitian vector bundle over a Riemann surface admits holomorphic line subbundles of arbitrarily large negative degree. 

\begin{lem}\label{lem:negative}
Let $\C^2 \to E \to \Sigma$ be a holomorphic vector bundle over a holomorphic curve.  For every line bundle $\mcl$ on $\Sigma$ of degree $k>0$ there exists $n \gg 1$ such that $\mcl^{-n} \subset E$ is a holomorphic subbundle of degree $-kn<0$. 
\end{lem}
\begin{proof}
Let $E_n=E \otimes \mcl^n$.  For $n$ large enough there exist holomorphic sections $(s_1, \ldots , s_r)$ of $E_n$ that span $(E_n)_{\sigma}$ for all $\sigma \in \Sigma$.  In this case there exists a nowhere vanishing section $s$ (see exercise 8.2 in \cite{hartshorne}) which defines an inclusion 
\be
0 \to \mco_{\Sigma} \to E_n
\ee
of the trivial line bundle $\mco_{\Sigma}$.  In turn this defines a holomorphic inclusion
\be
0 \to \mcl^{-n} \to E.
\ee
\end{proof}

The holomorphic line subbundle $\mcl \subset \mcw_+$ used to lift $\phi:\Sigma \to \s{6}$ to a pseudoholomorphic curve $\gamma: \Sigma \to \gro$ will be the pullback of the tautological bundle of $\gro$, $\mcl=\gamma^*(\mce)$.  Thus the degree of $\gamma$ is the degree of $\mcl$.  By Lemma \ref{lem:negative}, $\mcl$ can be chosen to have arbitrarily large negative degree. This shows that there is always a lift of $\phi:\Sigma \to \s{6}$ to a pseudoholomorphic curve with arbitrarily large negative degree, proving the first statement of Theorem \ref{thm:immersed}. 

Now we turn to the second statement.  Let $\phi:\Sigma \to \s{6}$ be a branched minimal immersion and $\mcl \subset \mcw_+$ a holomorphic subbundle.  Let $\gamma:=\mcl \circ \phi:\Sigma \to \gro$ be the holomorphic lift given by Lemma \ref{lem:hollift}.    Recall that $R_i=\{x \in \Sigma \;:\;{\rm I}_i(x)=0 \}$.  Then $\gamma$ is an immersion when 
\be\label{eqn:immersed}
R_1 \cap R_2 = \varnothing.
\ee
Locally one can adapt frames such that $\mcl=[f_0]$.  Since $\ffo$ only depends on $\theta_{od}$, it is an invariant of $\phi$ and so is independent of the choice of $\mcl$.  In contrast, $\ffe$ very much depends upon the geometry of $\mcl$.  For example, the second fundamental form of $\mcl$ is
\be
\ff_\mcl=\pi^{\perp}_{\mcl} \circ \nabla^{+} =f_1 \otimes f_0^* \otimes \frac{1}{2}\theta_{2}
\ee
and $\ffe$ also involves $\theta_2$.  To show that $R_1 \cap R_2 = \varnothing$ it is sufficient to find a holomorphic line subbundle $\mcl$ such that $\ff_{\mcl}$ does not vanish on the finite set of points $R_1$.   The existence of such an $\mcl$ is a consequence of the following lemma.

\begin{lem}\label{lem:immersion}
Let $\C^2 \to E \to \Sigma$ be a hermitian vector bundle on a holomorphic curve with connection $\nabla$ and hermitian inner product $(\;,\,)$.  Let 
\be
D=\{x'_1,\ldots,x'_m \}
\ee
be any finite set of points in $\Sigma$.  For any holomorphic line subbundle $L \subset E$ let 
\be
D_L=\{\sigma \in \Sigma : \ff_L(\sigma)=0 \}.
\ee  
Then there exists a holomorphic line subbundle $L \subset E$ such that $D \cap D_L=\varnothing$.
\end{lem}
\begin{proof}
Choose a meromorphic section $s$ of $E$ for which $s(x'_i) \neq 0,\infty$ for all $i=1 \ldots m$.  Let $V_i \subset E_{x'_i}$ be the subspaces spanned by $s(x'_i)$.  Let
\be
Z=\{t \in \Gamma(E): {\rm t\; is\; meromorphic}, \;t(x'_i)=0,\; \nabla t(x'_i) \pitchfork V_i \}\footnote{This set is always nonempty.  For instance let $t'$ be a meromorphic section of $E$ such that $\nabla t(x'_i) \pitchfork V_i$ for all $i=1 \ldots m$.  The set of such sections is open in the set of all meromorphic sections.  Choose a meromorphic function $u$ such that $u(x'_i)=0$ but $\ed u(x'_i) \neq 0$.  Then $t=t'u$ is in $Z$.}. 
\ee

Let $s_{\epsilon}=s+\epsilon t$ for any $t \in Z$, let $L_{\epsilon} \subset E$ be the line subbundle that it defines, let $\ff_{\epsilon}$ be its second fundamental form, and let $D_{\epsilon}=D_{L_{\epsilon}}$.  Define a map
\begin{align*}
M_{\epsilon}:&Z \to \oplus_{i=1}^{m} \left( E_{x'_i} \otimes T^*_{x'_i}\Sigma \right)
\end{align*}
where $M_{\epsilon}=(M_{\epsilon}^1,\ldots,M_{\epsilon}^m)$ and 
\be
M_{\epsilon}^i(t)=\ff_{L_{\epsilon}}(s_{\epsilon})(x'_i)
\ee
Due to the vanishing property of the $t$ we have
\be
M_{\epsilon}^i(t)=\ff_L(s)(x'_i)+\epsilon(\pi^{\perp}_L(\nabla t(x'_i)))
\ee
By definition $\pi^{\perp}_L(\nabla t(x'_i)) \neq 0$.  Thus for each $i=1 \ldots m$ there is a unique $\epsilon_i$ such that $M_{\epsilon_i}^i(t)=0$.  Choose $\epsilon' \neq \epsilon_i$ for all $i=1 \ldots m$.  Then $M^i_{\epsilon'}(t) \neq 0$ and thus $D_{\epsilon} \cap D = \varnothing$.
\end{proof}
\end{section}

\begin{section}{$2$-Ruled Cayley $4$-folds with a fixed asymptotic cone}\label{sec:deformations}
In this section we build upon the work of Bryant \cite{bryant2006}, Joyce \cite{joyce1r}, and Lotay \cite{lotayca} on ruled calibrated submanifolds.  We study the space of $2$-ruled Cayley $4$-folds $(M,\pi,\Sigma)$ that have a fixed asymptotic cone $M_0$.  This is a linear space, in fact the kernel of a  first-order linear differential operator \cite{lotayca}.  We will show that when the pseudoholomorphic curve defining a $2$-ruled Cayley cone $M_0$ has sufficiently negative degree, then $M_0$ is the asymptotic cone of some non-conical $2$-ruled Cayley $4$-fold.  Lotay obtained an analogous result in the case that the pseudoholomorphic curve is a $T^2$ of degree $0$.  If the pseudoholomorphic curve is embedded (so that the cone is smooth away from the origin) then this results in a smooth $2$-ruled Cayley $4$-fold $M$ with asymptotic cone $M_0$.  In fact, due to the linear structure, it gives rise to a one-paramter family of smooth $2$-ruled Cayley $4$-folds that degenerate to $M_0$ (See the paragraph after Lemma \ref{lem:lineareq}).  Finding a smooth Cayley $4$-fold $M$ with a fixed asymptotic cone is also one of the ingredients used in the process of desingularizing calibrated submanifolds with conical singularities \cite{slsing,lotaycodesing,hp}.     

Let 
\be
\gamma:\Sigma \to \gro
\ee
be a pseudoholomorphic curve and $U \subset \Sigma$ an open neighborhood on which one can write
\be
\gamma=e_1 \w e_2
\ee
for $e_1,e_2$ the first two legs of a $\Spin(7)$-adapted framing.  A $2$-ruled $4$-fold asymptotic to this cone can locally be parametrized as
\begin{align*}
\bx:&\R^2 \times U \to \Oc\\
&(r_1,r_2,\sigma) \mapsto r_1e_1+r_2e_2+s
\end{align*}   
where $s \in \ci{\mch}$.  This $4$-fold will be Cayley if 
\be
\bx^*(\psi_m)=0.
\ee
Let $s=e_a s^a$ for $8 \geq a,b,c \geq 3$ and also define the $1$-forms $s'^a:=\ed s^a+\om^a_b s^b.$  Using the fact that the asymptotic cone is already Cayley, the condition for $\bx$ to be Cayley simplifies to
\be
\begin{split}\label{eq:csec}
\psi_m(e_1,e_2,\ed s , \ed s)&=0\\
\psi_m(e_1,e_2,\ed e_i , \ed s)&=0
\end{split}
\ee
for $i=1,2$.  Notice that we may replace $\ed s$ by $\nabla^{\mch}s$ due to the congruence $\nabla^{\mch}s \equiv \ed s$ modulo $e_1,e_2$.

The solutions to these quadratic equations are the same as the solutions to a \emph{linear} $\db$-system which we now introduce.  Let 
\begin{align*}
a_1&=s^3+is^4\\
a_2&=s^6-is^7\\
a_3&=s^5-is^8
\end{align*}
so that if $\sigma=f\;a$ is a local section of $\mch'$ then $\Re(\sigma)=s$ is the corresponding section of $\mch$. Similarly, define the complex $1$-forms
\begin{align*}
\alpha_{1}&=s'^{3}+is'^{4}\\
\alpha_{2}&=s'^{6}-is'^{7}\\
\alpha_{3}&=s'^{5}-is'^{8}.
\end{align*}
Now recall that $\ip{\hspace{.1cm}}{\hspace{.01cm}}$ is the positive definite symmetric bilinear form on $\Oc$ (extended to $\Oc \otimes_{\R} \C$ as needed), so that $(v,w)=\ip{\ol{v}}{w}$ is a Hermitian inner product.
In terms of a local coframe we find that for $i=1,2,3$
\begin{align*}
\alpha_i&=\ip{\fb_i}{\nabla^{\mch}s}=\ip{\fb_i}{\ed \Re(\sigma)}\\
&=\Bigl \langle{\fb_i},{\Re \left( f(\kappa a+\ed a)-\frac{i}{2}\fb[\theta_{od}]a \right) }\Bigr \rangle \\
&=\ed a_i+ \kappa_{ij}a_j+\frac{i}{2}[\bar{\theta}_{od}]_{ij}\bar{a}_j.
\end{align*}

\begin{lem}\label{lem:lineareq}
Still using the notation above, the condition for $\bx=s+r_1e_1+r_2e_2$ to be Cayley \us{(}Equation \eqref{eq:csec}\us{)} is equivalent to the $\alpha_i$ being of type $(1,0)$ on $\Sigma$.
\end{lem}
\begin{proof}
Here is an outline.  Since $(\Sigma,\gamma)$ is pseudoholomorphic, there exist smooth functions $A_i:\Sigma \to \C$ such that $\theta_i=A_i\ed z$ where $z:\Sigma \to \C$ is a local holomorphic coordinate.   One must check that the equations in \eqref{eq:csec} are equivalent to the vanishing of the real and imaginary parts of $\alpha_i \w \ed z$, which is a routine calculation.   
\end{proof}
From this we see that the space of $2$-ruled Cayley $4$-folds with fixed asymptotic cone is the kernel of a first order linear differential operator.  Locally the equation is
\[
\db a_i+ \pi^{(0,1)}(\kappa_{ij})a_j+\frac{i}{2}[\bar{\theta}_{od}]_{ij}\bar{a}_j=0,
\]
where $\pi^{(0,1)}:\Om^{1} \to \Om^{(0,1)}$ is the projection.  So if $\bx=s+r_1e_1+r_2e_2$ is Cayley, then so is $\bx_{\lambda}=\lambda s +r_1e_1+r_2e_2$ for all $\lambda \in \R$.  Thus the limit of the image of $\bx_{\lambda}$ as $\lambda \to 0$ is $M_0$. It is in this sense that $\bx=\bx_{\lambda=1}$ is a desingularization of its asymptotic cone.

The fact that the Cayley equations reduce to a system of linear equations that are $\db$ to highest order makes one hope that they are equivalent to some functions being holomorphic.  This does not appear to be the case.  For instance, $\sigma$ is a holomorphic section of $\mch'$ if
\be
0=\db \sigma=\pi^{(0,1)}(\ed a +\kappa a ).
\ee
Unfortunately, the condition for $s$ to define a Cayley $4$-fold is not equivalent to this due to the $[\ol{\theta}_{od}]$ term in the expression for $\alpha_i$.  However, there is always a holomorphic line subbundle of $\mch'$ whose holomorphic sections do define Cayley $4$-folds.  

Let $\mcl \subset \mch'$ be the holomorphic line subbundle defined by $\ffo \in \Gamma(\mcv'_1)$ (see Section \ref{sec:fff}) and the equivalence ${\mathbb P}(\mch')={\mathbb P}(\mcv'_1)$.  In an adapted frame for which $\mcl=\C\cdot f_1$, we have $\theta_3=\theta_5=0$ and the simplification
\begin{align*}
\alpha_{1}&=\kappa_{11} a_1+\ed a_1\\
\alpha_{2}&=\kappa_{21} a_1\\
\alpha_{3}&=\kappa_{31} a_1,
\end{align*}
showing that the condition for such an $s \in \Gamma(\mcl)$ to define a Cayley $4$-fold is that $\sigma$ be a holomorphic section of $\mcl \subset \mch'$.

\begin{thm}
Fix a minimal immersion
\be
\phi:\Sigma \to \s{6}.
\ee
Then any pseudoholomorphic curve 
\be
\gamma:\Sigma \to \gro
\ee
with $\phi=\mcj \circ \gamma$ and sufficiently large negative degree always has a \emph{non-conical} $2$-ruled Cayley $4$-fold with asymptotic cone defined by $\gamma(\Sigma)$.
\end{thm}
\begin{proof}
We'll see that when $\deg(\gamma)$ is sufficiently large and negative then $\mcl$ admits many holomorphic sections.  

There is a decomposition
\be
\gamma^*(\Oc)=\gamma^*(\mce) \oplus \mcl \oplus \mcw.
\ee
where $\mcw=\mcl^{\perp} \subset \gamma^*(\mch)$.
Let 
\be
k=\gamma^*(c_1(\mcw')).
\ee
We have the relation
\be
c_1(\mcl)=-\deg(\gamma)-k.
\ee
As we choose different curves $\gamma$ that cover $\phi$, the integer $k$ remains unchanged.  Thus by requiring $\deg(\gamma)$ to be sufficiently negative, which we can achieve by Theorem \ref{thm:immersed}, Riemann-Roch implies that $h^0(\mcl)$ can be made as large as may be desired.  This guarantees holomorphic sections and thus deformations of the cone that are not conical.  
\end{proof}

\begin{rem}
Equation \eqref{eqn:immersed} and the proof of the Theorem \ref{thm:immersed} indicate that the generic pseudoholomorphic curve of large negative degree will be immersed.  One even expects it to be embedded since it is a real surface in a real $12$-dimensional space.  Thus we can expect that for $\deg(\gamma)$ sufficiently large and negative \us{(}and thus $h^0(\mcl) \gg 0$\us{)}, the generic holomorphic section will deform the cone into an \emph{immersed} $4$-fold.   In other words, we can expect that the generic $2$-ruled Cayley cone of large negative degree is the asymptotic cone of a smooth $2$-ruled Cayley $4$-fold.  This construction coincides with the one introduced by Lotay \cite{lotayca} in the case that $\gamma(\Sigma)$ is a torus and $\gamma^*(\mce)$ is trivial.   
\end{rem}
\end{section}

\appendix
\begin{section}{$\Spin(7)$ algebra}\label{sec:spin7}
Here are the details of the $\Spin(7)$ conventions used in this article.  The Cayley $4$-form used here is
\begin{equation}
\begin{split}
 \Phi=\,& \ed x^{5678}-\ed x^{5128} - \ed x^{5348}
-\ed x^{6138} -  \ed x^{6428} -\ed x^{7148} - \ed x^{7238}\\ +&\ed x^{1234}  - \ed x^{3467} - \ed x^{1267}
-\ed x^{2457} + \ed x^{1357} - \ed x^{2356} - \ed x^{1456}  
\end{split}
\end{equation}
where $\ed x^{ijkl}$ is shorthand for $\ed x^{i} \w \ed x^{j} \w \ed x^{k} \w \ed x^{l}$.  $\Phi$ is the real part of the quadruple cross product.  The imaginary components of the cross product
\be
\psi_m(e^0_i,e^0_j,e^0_k,e^0_l) = \ip{e^0_m}{e^0_i \times e^0_j \times e^0_k \times e^0_l}
\ee
for $1 \leq m \leq 7$ are 
\begin{align*}
&\psi_{1}=\ed x^{5137}+\ed x^{5864}+\ed x^{7123}+\ed x^{6421}+\ed x^{1675}+\ed x^{7538}+\ed x^{2786}+\ed x^{4283}\\
&\psi_{2}=\ed x^{3612}+\ed x^{2354}+\ed x^{4578}+\ed x^{7214}+\ed x^{5863}+\ed x^{3418}+\ed x^{5762}+\ed x^{1876}\\
&\psi_{3}=\ed x^{7314}+\ed x^{4821}+\ed x^{8476}+\ed x^{6342}+\ed x^{3756}+\ed x^{1253}+\ed x^{1578}+\ed x^{8256}\\
&\psi_{4}=\ed x^{1364}+\ed x^{6574}+\ed x^{4521}+\ed x^{3812}+\ed x^{2587}+\ed x^{6518}+\ed x^{7234}+\ed x^{6783}\\
&\psi_{5}=\ed x^{7352}+\ed x^{2847}+\ed x^{6814}+\ed x^{8632}+\ed x^{3817}+\ed x^{1653}+\ed x^{4571}+\ed x^{4265}\\
&\psi_{6}=\ed x^{5346}+\ed x^{6521}+\ed x^{5328}+\ed x^{4387}+\ed x^{1647}+\ed x^{7128}+\ed x^{2763}+\ed x^{4158}\\
&\psi_{7}=\ed x^{8463}+\ed x^{8531}+\ed x^{6274}+\ed x^{6137}+\ed x^{2485}+\ed x^{3574}+\ed x^{7521}+\ed x^{1628}
\end{align*}

The connection $\omega_{ij}$ is $\spin(7)$-valued if it satisfies
\be\label{eq:spin7rel}
\begin{split}
&\omega_{87}=-\omega_{65} +(\omega_{41} +\omega_{32})\\
&\omega_{86}=\omega_{75} +(\omega_{31} -\omega_{42})\\
&\omega_{85}=-\omega_{76}+( \omega_{43}+ \omega_{21})\\
&\omega_{81}=\omega_{74} +\omega_{63} +\omega_{52}\\
&\omega_{82}=\omega_{73} -\omega_{64} -\omega_{51}\\
&\omega_{83}=-\omega_{72}- \omega_{61}+ \omega_{54}\\
&\omega_{84}=-\omega_{71} +\omega_{62} -\omega_{53}.
\end{split}
\ee
\end{section}

\begin{section}{Holomorphic line subbundles over curves}\label{sec:sff}
Let $\mcw,\nabla$ be a hermitian vector bundle over a holomorphic curve $\Sigma$.  Any line subbundle $\mcl \subset \mcw$ has a second fundamental form whose definition we now recall.  Let $\pi:\Gamma(\mcw) \to \Gamma(\mcl)$ be the orthogonal projection and $\pi^{\perp}$ the complementary projection.
\begin{defn}
The second fundamental form of $\mcl \subset \mcw$ is a $\ci{\Sigma}$-linear bundle map
\[
\ff_{\mcl}:\Gamma(L) \to \Gamma({\mcl}^{\perp} \otimes T^*\Sigma).\\
\]
Given a section $s \in \Gamma(\mcl)$,
\[
\ff_{\mcl}(s)=\pi^{\perp} \circ \nabla (s).
\] 
\end{defn}
The main fact that will be used is contained in the following lemma.
\begin{lem}\label{lem:hlb}
$\mcl \subset \mcw$ is a \emph{holomorphic} line subbundle if and only if the image of $\ff_{\mcl}$ is contained in ${\mcl}^{\perp} \otimes \Omega^{(1,0)}(\Sigma)$.
\end{lem}
The criterion is equivalent to requiring that a holomorphic section of $\mcl \subset \mcw$ (with respect to the natural holomorphic structure induced by the hermitian structure $\mcl$ inherits from $\mcw$) is also a holomorphic section of $\mcw$ under the natural inclusion 
\[
\Gamma(\mcl) \hra \Gamma(\mcw).
\]
\end{section}
\bibliographystyle{amsplain}
\bibliography{2ruledcayleyconesbib}
\end{document}